\documentclass{article}
\usepackage{latexsym,amsmath,amssymb,graphicx}
\input amssym.def \input amssym

\newtheorem{te}{Theorem}[section]
\newtheorem{de}[te]{Definition}
\newtheorem{lm}[te]{Lemma}

\newtheorem{pp}[te]{Proposition}

\def\dokaz{\noindent{\bf Proof. }}
\def\kraj{\hfill $\Box$ \par \vspace*{2mm} }

\def\widemid{\hspace{1mm}\widetilde{\mid}\hspace{1mm}}
\def\nwidemid{\hspace{1mm}\widetilde{\nmid}\hspace{1mm}}
\def\rest{\upharpoonright}

\def\po{\exists}
\def\str{\rightarrow}
\def\Str{\Rightarrow}

\def\dl{\Leftrightarrow}
\def\ps{\subseteq}

\def\cl{{\rm cl}}
\def\cU{{\cal U}}
\def\cV{{\cal V}}

\begin{document}
\begin{center}
           {\large \bf Divisibility orders in $\beta N$}
\end{center}
\begin{center}
{\small \bf Boris  \v Sobot}\\[2mm]
{\small  Department of Mathematics and Informatics, University of Novi Sad,\\
Trg Dositeja Obradovi\'ca 4, 21000 Novi Sad, Serbia\\
e-mail: sobot@dmi.uns.ac.rs}
\end{center}
\begin{abstract} \noindent
We consider five divisibility orders on the Stone-\v Cech compactification $\beta N$. We find some possible lengths of chains and antichains and number of maximal and minimal elements, as well as some other ordering properties of these relations.
\vspace{1mm}\\


\end{abstract}

\section{Introduction}

The book \cite{HS} considers extensions of semigroup operations on discrete spaces $S$ to their Stone-\v Cech compatifications $\beta S$. We are interested in the extension of the multiplication on the set $N$ of natural numbers in this way. The operation $\cdot$ on $N$ is extended to $\beta N$ as follows:
$$A\in p\cdot q\dl\{n\in N:A/n\in q\}\in p,$$
where $A/n=\{\frac an:a\in A,n\mid a\}$. In particular, if $n\in N$ and $q\in\beta N$ then $A\in nq$ if and only if $q\in A/n$. The topology on $\beta N$ is defined by taking ${\bar A}=\{p\in\beta N:A\in p\}$ (for $A\subseteq N$) as base sets (the set ${\bar A}$ is the closure of $A$ so there is no abuse of notation).\\

Let us fix some notation. Identifying elements of $N$ with the corresponding principal ultrafilters, we will denote $N^*=\beta N\setminus N$. The (unique) continuous extension of a function $f:N\str N$ to $\beta N$ will be denoted by $\widetilde{f}$. $\mid$ is the divisibility relation on $N$, and $\mid[A]=\{m\in N:\po a\in A\;a\mid m\}$. Let also $\cU=\{S\ps N:S\mbox{ is upward closed for }\mid\}$ and $\cV=\{S\ps N:S\mbox{ is downward closed for }\mid\}$. We also mention that, since almost all the results we use are contained in the book \cite{HS}, for readers' convenience we chose to cite the book instead of various papers in which results may have appeared first.\\

In \cite{S1} we defined four possible extensions of the divisibility relation $\mid$ on $N$ to $\beta N$. The eventual goal of investigating these relations is to try to translate problems from elementary number theory (of infinite character, i.e.\ problems dealing with infinity of certain subsets of $N$) into $\beta N$ and use topological methods to approach them.

\begin{de}
Let $p,q\in\beta N$.

(a) $q$ is left-divisible by $p$, $p\mid_L q$, if there is $r\in\beta N$ such that $q=rp$.

(b) $q$ is right-divisible by $p$, $p\mid_R q$, if there is $r\in\beta N$ such that $q=pr$.

(c) $q$ is mid-divisible by $p$, $p\mid_M q$, if there are $r,s\in\beta N$ such that $q=rps$.

(d) $p\widemid q$ if for all $A\in p$, $\mid[A]\in q$ holds.
\end{de}

In semigroup theory $\mid_L$ and $\mid_R$ are known as the Green relations: they are equivalent to the inclusion relation on the sets of principal left (or right) ideals, for example $p\mid_L q$ if and only if $\beta Nq\subseteq\beta Np$. Hence they have been considered before. The relation $\widemid$ was introduced by analogy with functions $\widetilde{f}$, extending $\mid$ in such way to satisfy certain continuity conditions, and in \cite{S1} it was proved that $\widemid$ is the maximal extension of $\mid$ which is continuous in that sense. In this paper we will investigate some ordering properties of these relations, adding one more, $\mid_{LN}$.\\

All the relations $\mid_L$, $\mid_R$, $\mid_M$ and $\widemid$ are preorders (reflexive and transitive), but none is antisymmetric (see Section 4 of \cite{S1}). So for each of them we introduce another relation: $p=_Lq$ if $p\mid_Lq$ and $q\mid_Lp$, and $=_R$, $=_M$ and $=_{\sim}$ are defined analogously. All these are equivalence relations, and all the divisibility relations can be viewed as partial orders on respective factor sets (we will use the same notation for orders on factor sets as for preorders above). Respective equivalence classes are denoted by $[p]_L$, $[p]_R$, $[p]_M$ and $[p]_\sim$.

\begin{lm}\label{cancel}
(a) If $p$ is right cancelable then $[p]_L=\{p\}$. (b) If $p$ is left cancelable then $[p]_R=\{p\}$.
\end{lm}

\dokaz (a) Assume the opposite, that there is $q\neq p$ such that $p=_Lq$. This means that $p=xq$ and $q=yp$ for some $x,y\in\beta N$. Then $p=xyp$, so since $p$ is right cancelable, $xy=1$. But $N^*$ is an ideal of $\beta N$ (\cite{HS}, Theorem 4.36), which means that $x=y=1$, so $p=q$.

(b) is proven analogously.\kraj

Note also that the sets of right cancelable and left cancelable elements are downward closed in $\mid_R$ and $\mid_L$ respectively.\\

In the next proposition we collect several useful facts concerning elements of $N$.

\begin{pp}\label{deljN}
(a) (\cite{HS}, Theorem 6.10) Elements of $N$ commute with all elements of $\beta N$.


(b) (\cite{HS}, Lemma 6.28) If $m,n\in N$ and $p\in\beta N$, then $mp=np$ implies $m=n$.

(c) (\cite{S1}, Lemma 5.1) If $n\in N$, each of the statements: (i) $n\mid_L p$, (ii) $n\mid_R p$, (iii) $n\mid_M p$, (iv) $n\widemid p$ and (v) $nN\in p$ are equivalent.
\end{pp}

The previous lemma allows us to drop subscripts and write only $n\mid p$ for $n\in N$.

\begin{lm}\label{Dodp}
The following conditions are equivalent:

(i) $p\widemid q$;

(ii) $D(p)\ps q$;

(iii) $p\cap\cU\ps q\cap\cU$;

(iv) $q\cap\cV\ps p\cap\cV$.

\end{lm}

\dokaz The equivalence of (i), (ii) and (iii) was proved in \cite{S1}, Theorem 6.2.

(iii)$\dl$(iv) follows easily from the fact that $A\in\cU$ iff $A^c\in\cV$.\kraj

\begin{lm}
(a) For each $p\in\beta N$ the sets $p\uparrow_L=\{q:p\mid_L q\}$ are closed;

(b) for each $p\in\beta N$ the sets $p\uparrow_\sim=\{q:p\widemid q\}$, $p\downarrow_\sim=\{q:q\widemid p\}$ and $[p]_\sim$ are closed.
\end{lm}

\dokaz (a) $p\uparrow_L=\beta Np=\overline{Np}$, which is clearly closed.

(b) By Lemma \ref{Dodp}(ii) $q\in p\uparrow_\sim$ iff $q\in\bigcap_{A\in D(p)}\overline{A}$, which is a closed set. By Lemma \ref{Dodp}(iv) $p\downarrow_\sim=\bigcap_{B\in p\cap\cV}\overline{B}$, also a closed set. $[p]_\sim=p\uparrow_\sim\cap p\downarrow_\sim$, so it is also closed.\kraj

In \cite{HS}, Definition 1.34, orders $\geq_L$, $\geq_R$ and $\geq$ were defined on the set $E(\beta N)$ of idempotents of $(\beta N,\cdot)$: $p\geq_L q$ if $pq=q$, $p\geq_R q$ if $qp=q$, and $p\geq q$ if both $p\geq_L q$ and $p\geq_R q$.

\begin{lm}
$\geq_R$ is the restriction of $\mid_L$ to $E(\beta N)$, and $\geq_L$ is the restriction of $\mid_R$ to $E(\beta N)$.
\end{lm}

\dokaz We prove the result for $\mid_L$; the proof for $\mid_R$ is analogous. Clearly, for $p,q\in E(\beta N)$ $p\geq_R q$ implies $p\mid_L q$. On the other hand, if $p\mid_L q$ then there is $x\in\beta N$ such that $q=xp$. Then $qp=xpp=xp=q$ because $p$ is an idempotent.\kraj

Clearly, if for two of the considered relations $\tau$ and $\sigma$ holds $\tau\ps\sigma$, then equivalence classes of $=_\sigma$ are unions of equivalence classes of $=_\tau$. We will show that the following diagram holds:

$$\begin{array}{c}
\\
\mid_L\\
\\
\mid_R
\end{array}\hspace{-3mm}
\begin{array}{c}
\rotatebox[origin=c]{45}{$\subset$}\\
\rotatebox[origin=c]{-45}{$\subset$}\\
\rotatebox[origin=c]{45}{$\subset$}
\end{array}\hspace{-2mm}
\begin{array}{l}
\mid_{LN}\\
\\
\mid_M\;\subset\;\widemid\\
\\
\end{array}
$$

The part of the diagram concerning $\mid_L$, $\mid_R$, $\mid_M$ and $\widemid$ was explained in \cite{S1}. The inclusion $\mid_L\subset\mid_{LN}$ will be clear from the definition of $\mid_{LN}$. Why $\mid_R$, $\mid_M$ and $\widemid$ are incomparable with $\mid_{LN}$ will be explained at the end of Section \ref{minmax}.

\section{The preorder $\mid_{LN}$}

\begin{pp}\label{mozeN}
(a) (\cite{HS}, Theorem 3.40) If $A$ and $B$ are countable subsets of $\beta N$ such that ${\bar A}\cap{\bar B}\neq\emptyset$, then $A\cap{\bar B}\neq\emptyset$ or ${\bar A}\cap B\neq\emptyset$.

(b) (\cite{HS}, Theorem 6.19) If $xp=yq$ for $p,q,x,y\in\beta N$, then there are $n\in N$ and $z\in\beta N$ such that either $np=zq$ or $zp=nq$.
\end{pp}

The previous result suggests introduction of another relation on $\beta N$, representing divisibility "up to elements of $N$".

\begin{de}\label{defLN}
$p\mid_{LN} q$ if there is $n\in N$ such that $p\mid_Lnq$.
\end{de}

To see that this relation is strictly stronger than $\mid_L$, note that for any $p\in N^*$ not divisible by 2 we have $2p\mid_{LN} p$, but not $2p\mid_L p$.\\

The relation $\mid_{LN}$ is also a preorder so we introduce $=_{LN}$ and $[p]_{LN}$ as for other relations.

\begin{lm}
For each $p\in \beta N$, $[p]_{LN}=\{mr:m\in N,r\in\beta N,\exists n\in N(r\mid_L p\land p\mid_L nr)\}$.
\end{lm}

\dokaz It is obvious that, if there is $n\in N$ such that $r\mid_L p$ and $p\mid_L nr$, then $mr=_{LN}p$. So let $p\mid_{LN}q$ and $q\mid_{LN}p$ for some $q\in\beta N$. $q\mid_{LN}p$ means that there is $m\in N$ such that $q\mid_L mp$, i.e.\ $xq=mp$ for some $x\in\beta N$. We can assume that $m$ is minimal such element of $N$; let us show that this implies $m\mid q$.

If not, there is a prime $k\mid m$ such that $k\nmid q$, so since prime numbers are also prime in $\beta N$ (\cite{S1}, Lemma 7.3), and $k\mid xq$, we would have $k\mid x$, i.e.\ $x=kx_1$. Then we could cancel out $k$ and get $x_1q=\frac mkp$, which is a contradiction with the minimality of $m$.

So $q=mr$ for some $r\in\beta N$, and $r\mid_L p$. That $p\mid_L nr$ for some $n\in N$ follows directly from $p\mid_{LN}q$.\kraj

The following lemma is what may make this new relation useful.

\begin{lm}\label{linLN}
For every $q\in\beta N$ the set $\{[p]_{LN}:p\mid_{LN}q\}$ is linearly ordered by $\mid_{LN}$.
\end{lm}

\dokaz Let $p_1\mid_{LN}q$ and $p_2\mid_{LN}q$. This means that there are $n_1,n_2\in N$ and $x_1,x_2\in\beta N$ such that $n_1q=x_1p_1$ and $n_2q=x_2p_2$. Then $n_2x_1p_1=n_1n_2q=n_1x_2p_2$. By Proposition \ref{mozeN}(b) one of the elements $p_1$ and $p_2$ must be $\mid_{LN}$-divisible by the other.\kraj

\section{Chains and antichains}

In this section we investigate possible lengths of chains and antichains in our partial orders.

\begin{pp}
(\cite{HS}, Theorem 6.73) There is an infinite strictly $\mid_R$-descending chain in $\beta N$.
\end{pp}

\begin{lm}
There is an infinite strictly $\mid_L$-descending chain of right cancelable elements in $\beta N$.
\end{lm}

\dokaz We construct the wanted chain as follows: let $p_0\in\bigcap_{n\in N}\overline{2^nN}$ be right cancelable (by \cite{S1}, Theorem 5.2, this set is a nonempty $G_\delta$ set, thus it contains an open subset, and by \cite{HS}, Theorem 8.10 it contains a right cancelable element). Let $p_n=2p_{n+1}$ for $n\in\omega$. By induction it is easy to prove that all $p_n$ are right cancelable: for example $xp_1=yp_1\Str xp_0=yp_0\Str x=y$. These elements are also different (by Proposition \ref{deljN}(b)), belong to different $\mid_L$-classes and $\dots p_2\mid_L p_1\mid_L p_0$.\kraj


\begin{lm}
There is an strictly $\widemid$-descending chain of length $\omega+1$.
\end{lm}

\dokaz Let $P$ be the set of prime numbers, and let $\langle P_n:n<\omega\rangle$ be a sequence of sets such that $P_0=P$, $P_{n+1}\subset P_n$ and $P_n\setminus P_{n+1}$ is infinite for all $n\in\omega$. For each $n<\omega$ let $X_n=\{k\in N:\mbox{all prime divisors of }k\mbox{ belong to }P_n\}$.

For each $n<\omega$ the set $\{X_n\}\cup\{A\in\cU:A\cap X_n\neq\emptyset\}$ has the finite intersection property: let $A_1,A_2,\dots,A_k\in\cU$ be given with nonempty intersections with $X_n$. For every $i=1,2,\dots,k$ choose an element $a_i\in A_i\cap X_n$; then $LCM(a_1,a_2,\dots,a_k)\in A_1\cap A_2\cap\dots\cap A_k\cap X_n$. Note that $A\cap X_n\neq\emptyset$ for $A\in\cU$ actually means that $A\cap X_n$ is infinite.

Hence we can pick ultrafilters $p_n$ so that $p_n\cap\cV=\{A\in\cV:P_n\subseteq A\}$. Hence $p_m\cap\cV\subset p_n\cap\cV$ for $m<n<\omega$ so, by Lemma \ref{Dodp}, $p_n\widemid p_m$ and $p_m\nwidemid p_n$.

Finally, the family $\bigcup_{n<\omega}(p_n\cap\cV)$ has the finite intersection property, so there is an ultrafilter containing $\bigcup_{n<\omega}(p_n\cap\cV)$, which is below all $p_n$ for $n<\omega$.\kraj

In the next two lemmas we adapt ideas from the proof of Lemma 9.22 of \cite{HS}.

\begin{lm}\label{ascendL}
Every $\mid_L$-ascending chain of length $\omega$ has an upper bound in $\beta N$.
\end{lm}

\dokaz Let $\langle q_n:n\in\omega\rangle$ be a $\mid_L$-ascending chain. Let $q\in\cl\{q_n:n\in\omega\}\setminus\{q_n:n\in\omega\}$ be arbitrary. Then, for each $m\in\omega$, $q_n\in\beta N q_m$ for all $n\geq m$ so, since $\beta Nq_m$ is closed, $q\in\cl\{q_n:n\geq m\}\subseteq\beta Nq_m$ i.e.\ $q_m\mid_L q$.\kraj

\begin{lm}
Every strictly $\mid_{LN}$-ascending chain of length $\omega$ has an upper bound $q$ in $\beta N$ that is right cancelable.
\end{lm}

\dokaz Let $\langle r_n:n\in\omega\rangle$ be a strictly $\mid_{LN}$-ascending sequence, and for each $n\in\omega$ let $k_n\in N$ be such that $r_i\mid_L k_nr_n$ for all $i<n$. Let $q_n=k_nr_n$. Then $\langle q_n:n\in\omega\rangle$ is a $\mid_L$-ascending sequence which is also strictly $\mid_{LN}$-ascending. As in Lemma \ref{ascendL} we can find $q\in\cl\{q_n:n\in\omega\}$ which is a $\mid_L$-upper bound of $\langle q_n:n\in\omega\rangle$, and hence a $\mid_{LN}$-upper bound of $\langle r_n:n\in\omega\rangle$ as well.

Suppose $q$ is not right cancelable. Then, by Theorem 8.11 of \cite{HS}, there are $x\in N^*$ and $a\in N$ such that $xq=aq$. Hence, $aq\in\cl\{aq_n:n\in\omega\}\cap\cl((N\setminus\{a\})q)$ so, by Proposition \ref{mozeN}(a), we have one of the following two possibilities: either $\{aq_n:n\in\omega\}\cap\cl((N\setminus\{a\})q)\neq\emptyset$ or $\cl\{aq_n:n\in\omega\}\cap (N\setminus\{a\})q\neq\emptyset$. The first one leads to contradiction right away, since $aq_n=yq$ would mean that $q\mid_{LN} q_n$ and thus $q_{n+1}\mid_{LN} q_n$ as well.

So $aq'=bq$ for some $q'\in\cl\{q_n:n\in\omega\}$ and some $b\in N\setminus\{a\}$. This means that $\cl\{aq_n:n\in\omega\}\cap\cl\{bq_n:n\in\omega\}\neq\emptyset$ so, again by Proposition \ref{mozeN}(a), either $aq_m=bq''$ or $aq''=bq_m$ for some $m\in N$ and some $q''\in\cl\{q_n:n\in\omega\}$. ($q''\notin\{q_n:n<\omega\}$ because $q''=q_n$ would imply that $q_m=_{LN}q_n$.) Without loss of generality assume the first possibility. It follows that $q''\mid_{LN} q_m$, so $q_{m+1}\mid_{LN}q_m$, a contradiction again.\kraj

For divisibility relations we will consider the following notion of incompatibility.

\begin{de}
If $\rho$ is a preorder on $\beta N$ we will say that two elements $x$ and $y$ are $\rho$-compatible if there is $z\neq 1$ such that $z\rho x$ and $z\rho y$.

We define the compatibility relation $C_L$ on $N^*$ as follows: $pC_L q$ if there is $r\in N^*$ such that $r\mid_L p$ and $r\mid_L q$. The relation $\sim_L$ is its transitive closure: $p\sim_L q$ if there are $x_1,x_2,\dots,x_k\in\beta N$ such that $pC_L x_1,x_1C_L x_2,\dots,x_kC_L q$.
\end{de}

Clearly, $p\mid_L q$ implies $pC_L q$ (and hence $p\sim_L q$). $\sim_L$ is an equivalence relation.\\

In \cite{HS}, Definition 6.48, a relation $R$ on $N^*$ is defined by $pRq\dl\beta Np\cap\beta Nq\neq\emptyset$; hence $pRq$ if and only if $p$ and $q$ are $\mid_L^{-1}$-compatible. Then $\sim_L$ is also the transitive closure of $R$. The graphs of relations $\mid_L$, $\mid_{LN}$ (more precisely, of their symmetric closures), $R$, $C_L$ and $\sim_L$ on $N^*$ have the same connected components. The following result (a straightforward consequence of Theorem 6.53 of \cite{HS} and the considerations above) shows that these graphs are far from being connected.

\begin{pp}
There are $2^{\goth c}$ equivalence classes $[p]_{\sim_L}=\{q\in N^*:p\sim_L q\}$. Each of these classes is nowhere dense in $N^*$ and it is a left ideal of $\beta N$.
\end{pp}

For the relation $\widemid$ (and consequently for $\mid_R$ and $\mid_M$) we have a weaker result. We remind the reader that sets $A$ and $B$ are almost disjoint if $A\cap B$ is finite, and that on any infinite set there exists an almost disjoint family of cardinality $\goth c$ (see, for example, Lemma 3.1.2 of \cite{vM}).

\begin{te}
There is a family $\{p_\alpha:\alpha<{\goth c}\}$ of $\widemid$-incompatible elements in $\beta N$.
\end{te}

\dokaz Let $\{A_\alpha:\alpha<{\goth c}\}$ be an almost disjoint family of infinite subsets of the set $P$ of prime numbers. For $\alpha<{\goth c}$ let $p_\alpha\in N^*$ be an ultrafilter containing $A_\alpha$. We will prove that any two ultrafilters $p_\alpha$ and $p_\beta$ for $\alpha\neq\beta$ are $\widemid$-incompatible. First, since each of the sets $nN$ for $n\in N\setminus\{1\}$ is almost disjoint with $A_\alpha$, it follows that $nN\notin p_\alpha$ so by Proposition \ref{deljN}(c) $p_\alpha$ is not divisible by any $n\in N\setminus\{1\}$.

Now assume there is $r\in N^*$ such that $r\widemid p_\alpha$ and $r\widemid p_\beta$. Let $B_\alpha=\;\mid[P\setminus A_\alpha]$ and $B_\beta=\;\mid[P\setminus A_\beta]$; then $B_\alpha,B_\beta\in\cU$. Since $B_\alpha$ is disjoint from $A_\alpha$, $B_\alpha\notin p_\alpha$ so by Lemma \ref{Dodp} $B_\alpha\notin r$ as well. In the same way we conclude $B_\beta\notin r$. Also, $S=N\setminus(P\cup\{1\})$ is a set in $\cU$ disjoint from $A_\alpha$, so again $S\notin r$ and $P\in r$. This means that $P\setminus(B_\alpha\cup B_\beta)=A_\alpha\cap A_\beta$ must be in $r$, so $r$ must be a principal ultrafilter; a contradiction.\kraj

\section{Maximal and minimal elements}\label{minmax}

In \cite{HS}, Theorems 1.51 and 1.64, it is shown that $\beta N$ has the smallest ideal, denoted by $K(\beta N)$ and that $K(\beta N)=\bigcup\{L:L\mbox{ is a minimal left ideal of }\beta N\}=\bigcup\{R:R\mbox{ is a minimal right ideal of }\beta N\}$. Clearly, every minimal left ideal $L$ is principal and moreover generated by any element $p\in L$.

\begin{te}\label{max}
(a) $\mid_L$ has $2^{\goth c}$ maximal classes, they are exactly minimal left ideals of $\beta N$, and for every $q\in\beta N$ there is $p$ such that $q\mid_L p$ and $[p]_L$ is maximal.

(b) $\mid_R$ has $2^{\goth c}$ maximal classes, they are exactly minimal right ideals of $\beta N$, and for every $q\in\beta N$ there is $p$ such that $q\mid_R p$ and $[p]_R$ is maximal.

(c) For every $q\in\beta N$ there is $p$ such that $q\mid_{LN} p$ and $[p]_{LN}$ is maximal.

(d) $(\beta N/_{=_M},\mid_M)$ has the greatest element, and it is exactly the class $K(\beta N)$.

(e) $(\beta N/_{=_\sim},\widemid)$ has the greatest element, $\bigcap_{A\in\cU}\overline{A}$, containing $K(\beta N)$.
\end{te}

\dokaz (a) By Theorem 6.44 of \cite{HS} there are $2^{\goth c}$ minimal left ideals (and each of them has exactly $2^{\goth c}$ elements). But $[p]_L$ is maximal if and only if $\beta Np$ is a minimal left ideal. Since every principal left ideal contains a minimal left ideal (\cite{HS}, Corollary 2.6), there is a $\mid_L$-maximal element above every $q\in\beta N$.

(b) is proved analogously to (a), using Corollary 6.41 of \cite{HS}.

(c) If for some $p,q\in\beta N$ we have $p\mid_{LN}q$ and $p\nmid_L q$, then there is $r=_{LN}q$ such that $p\mid_L r$. We conclude that $\mid_L$-maximal elements are also $\mid_{LN}$-maximal so, by (a), above every $q\in\beta N$ there is a $\mid_{LN}$-maximal element.


(d) We first prove that all elements of $K(\beta N)$ are in the same $=_M$-equivalence class. Let $p,q\in K(\beta N)$. By Theorem 2.7(d) of \cite{HS} the left ideal $\beta Np$ intersects the right ideal $q\beta N$; let $r\in \beta Np\cap q\beta N$. Then $p=_L r$ (because $\beta Np=\beta Nr$) and $r=_R q$, so $p=_M r=_M q$.

It remains to prove that no element $p\notin K(\beta N)$ is in this maximal class (or above it). Assume the opposite, that there is $p$ such that $q\mid_M p$ for $q\in K(\beta N)$. But this means that $p=aqb$ for some $a,b\in\beta N$, and since $K(\beta N)$ is an ideal, it follows that $p\in K(\beta N)$ as well.

(e) Since $\cU$ has the finite intersection property, there are ultrafilters containing all the sets from $\cU$, and by Lemma \ref{Dodp} they are clearly maximal for $\widemid$. Since $\mid_M\subseteq\widemid$, by (d) all the ultrafilters from $K(\beta N)$ are among them.\kraj

\begin{lm}
If $p\in\beta N$ is right cancelable, then it has $2^{\goth c}$ incomparable $\mid_L$-successors (and $2^{\goth c}$ incomparable $\mid_R$-successors).
\end{lm}

\dokaz Let $q_\alpha$ (for $\alpha<2^{\goth c}$) be $\mid_L$-maximal elements such that $q_\alpha\beta N$ are different minimal left ideals. Then $q_\alpha p\nmid_L q_\gamma p$ for $\alpha\neq\gamma$: if $q_\gamma p=xq_\alpha p$, by cancelability we would have $q_\gamma=xq_\alpha$, so $q_\gamma\in\beta Nq_\alpha$, and by minimality $\beta Nq_\alpha$ and $\beta Nq_\gamma$ would be the same minimal left ideals.\kraj

Each of the orders we are investigating clearly has the smallest element: the one-element class $[1]=\{1\}$ is the smallest in $\mid_L$, $\mid_R$, $\mid_M$ and $\widemid$, and the equivalence class $[1]=N$ is the smallest in the order $\mid_{LN}$. It is more interesting to ignore the class $[1]$ and define, for each relation $\rho$, the set of minimal elements $M_\rho$ to be the union of minimal classes in $(\beta N/_{=_\rho}\setminus [1],\rho)$.\\

\begin{lm}
(a) $M_{\mid_L}=M_{\mid_R}=M_{\mid_M}$: it is the set of irreducible ultrafilters (those that can not be written as $pq$ for $p,q\in\beta N\setminus\{1\}$).

(b) For $\rho\in\{\mid_L,\mid_R,\mid_M,\widemid\}$, $M_\rho\cap N=P$.

(c) For $\rho\in\{\mid_L,\mid_R,\mid_M\}$, $P^*\subset M_\rho\cap N^*$.

(d) There are $\mid_{LN}$-minimal elements.
\end{lm}

\dokaz (a) is obvious.

(b) is obvious for $\rho\in\{\mid_L,\mid_R,\mid_M\}$. But $\widemid\rest_{N^2}=\mid_L\rest_{N^2}$, so the result holds for $\widemid$ too.

(c) $P^*\subseteq M_\rho\cap N^*$ follows from Theorem 7.3 and the strict inclusion from Theorem 7.5 of \cite{S1}.

(d) follows from Theorem 8.22 of \cite{HS}.\kraj

Now $\mid_M\not\ps\mid_{LN}$ is clear since by Theorem \ref{max} the order $\mid_M$ has the greatest element and $\mid_{LN}$ does not. But then $\mid_R\not\ps\mid_{LN}$ as well, since $\mid_R\ps\mid_{LN}$ along with $\mid_L\ps\mid_{LN}$ would imply that the transitive closure of $\mid_L\cup\mid_R$ (which is $\mid_M$) would be contained in $\mid_{LN}$ too.

Finally, $\mid_{LN}\not\ps\widemid$ because $2\mid_{LN}1$ but not $2\widemid 1$.

\end{document}